\documentclass[12pt]{article}
\usepackage{amssymb,amsmath,epsfig}

\newtheorem{theorem}{Theorem}[section]
\newtheorem{thm}[theorem]{Theorem}

\newtheorem{prop}[theorem]{Proposition}

\newtheorem{lemma}[theorem]{Lemma}
\newtheorem{cor}[theorem]{Corollary}

\newtheorem{definition}{Definition}[section]

\newtheorem{remark}{Remark}[section]

\newenvironment{proof}{\par\medskip\noindent{\em Proof. }}{\hfill $\square$\par\medskip}

\def\<{\langle}
\def\>{\rangle}
\def\Z{\mathbb{Z}}
\def\P{\mathcal P}
\def\Q{\mathbb{Q}}
\def\N{\mathbb{N}}
\def\T{\mathbb{T}}
\def\G{\Gamma}
\def\g{\gamma}
\def\o{\omega}
\def\a{\alpha}
\def\b{\beta}
\def\l{\lambda}
\def\FP{{\rm{FP}}}
\def\ssm{\smallsetminus}

\title{Subgroups of direct products of elementarily free groups}
\author{Martin R. Bridson and James Howie}
\begin{document}
\date{}
\maketitle

\begin{abstract} We exploit Zlil Sela's description of 
the structure of groups having the same 
elementary theory as free groups:
they and their
finitely generated subgroups form a prescribed
subclass $\mathcal E$ of the hyperbolic limit groups.
We prove that if $G_1,\dots,G_n$ are in $\mathcal E$
then a subgroup $\Gamma\subset G_1\times\dots\times G_n$
is of type $\FP_n$ if and only if $\Gamma$ is
itself, up to finite index, 
the direct product of at most 
$n$ groups from $\mathcal E$. This
answers a question of Sela.
\end{abstract}
Examples of Stallings \cite{stall} and Bieri \cite{bieri1}
show that finitely presented subgroups of a direct product of finitely
many free groups can be rather wild. In contrast, Baumslag and Roseblade
\cite{BR} proved that the only  finitely presented subgroups of a direct
product of  {\em{two}} free groups   are the obvious ones: such
a subgroup is either free or else has a subgroup of finite index that
is the product of its intersections with the factors. 

In 
\cite{BHMS} Miller, Short and the present authors explained
this apparent contrast 
by showing that all of the exotic behaviour among subdirect
products of free groups arises from a lack of homological finiteness. More
precisely, if a subgroup $S$ of a direct product of $n$ free groups has finitely
generated homology up to dimension
$n$, then $S$ has a subgroup of finite index that is
isomorphic to a direct product of free groups.

We proved a similar theorem for subdirect products
of surface groups. This has implications for  the
understanding of the fundamental groups of compact K\"ahler manifolds. Indeed
the remarkable work of Delzant and Gromov \cite{DG} shows that if the
fundamental group $\G$ of such a  manifold  is torsion-free and
has sufficiently many multi-ended
splittings, then there is a short exact sequence $1\to \Z^n\to\G_0\to S\to 1$,
where $S$ is a subdirect product of surface groups and $\G_0$ is a 
subgroup of finite index in $\G$.

Different
attributes of surface groups enter the proof in \cite{BHMS} in subtle ways
making it difficult to assess which are vital and which are artifacts of the proof.
Examples show that the splitting phenomenon for 
subdirect products does not extend to arbitrary
2-dimensional hyperbolic groups (even small cancellation
groups), nor to  fundamental groups of higher-dimensional hyperbolic
manifolds \cite{mrb}. However, advances in geometric group theory during the
last few years suggest that there is a class of groups which is 
much more fundamentally tied to surface groups than either of the
these classes, namely {\em{limit groups}}.

Kharlampovich and Myasnikov have studied limit groups extensively
under the name {\em{fully residually free groups}}
\cite{KM1}, \cite{KM2}.  Remeslennikov \cite{rem} had previously
referred to them as $\exists$-free groups, reflecting
the fact that these are precisely the groups that have the same 
existential theory as a free group. The name {\em{limit
group}} was introduced by Zlil Sela to emphasise the fact that these are
precisely the class of groups that arise when one takes limits of 
stable sequences of homomorphisms $\phi_n:G\to F$, where $G$ is an arbitrary
finitely generated group and $F$ is a free group (see section 1).

The existential theory of a group $G$
is the set of first order sentences in the language of group
theory that contain only one quantifier $\exists$ and are true in $G$.
 The {\em{elementary
theory}} of $G$ is the set of all first order sentences that are
true in $G$. Famously, Alfred
Tarski asked which groups had the same elementary theory as a
free group. This problem was solved by Zlil Sela  (\cite{S1} to \cite{S6}).
At about the same time,
a somewhat different solution was proposed by Kharlampovich and Myasnikov,
but our results here rely on Sela's treatment, in particular his hierarchical
description of this class of groups (Definition \ref{d:sela}). 
We refer to a group as {\em{elementarily free}} if it has the same elementary theory
as a free group.  This is an important subclass of the class of limit groups.

The following is our main result.

\begin{theorem}\label{t:main}
Let $G_1,\dots,G_n$ be (subgroups of)
elementarily free groups and let
 $\G\subset G_1\times\cdots\times G_n$ be a subgroup.
Then
$\G$ is of type
$\FP_n(\mathbb Q)$ if and only if  there are finitely generated
subgroups $H_i\subset G_i$ for $i=1,\dots,n$
such that $\G$ is isomorphic to
a  subgroup of finite-index in $H_1\times\cdots\times H_n$.
\end{theorem}

For simplicity of exposition, all homology groups considered in this
paper will be with trivial coefficient module $\Q$, so it is natural
to use the finiteness condition $\FP_n(\Q)$ in the statement of our theorem.
However, with minor modifications, our arguments also apply
with other trivial
coefficient modules, giving corresponding results under
the finiteness conditions $\FP_n(R)$ for other
suitable rings $R$.

Theorem  \ref{t:main} answers a 
question of Sela \cite[I(12)]{S9}, who also asked if the
stronger version of the above result holds, 
in which `elementarily free groups'
is replaced by `limit groups'.  We conjecture that 
the stronger result is indeed
true, but it cannot be proved using the methods of this paper alone.
In a forthcoming article \cite{BHo2}, 
we use a
different technique to prove the
stronger conjecture in the case $n=2$.  

By combining a slight refinement of Theorem \ref{t:main} (see \ref{last})
 with Theorem 3 of
\cite{BHo1} we obtain:

\begin{theorem}\label{intersects}
Let $G_1,\dots,G_n$ be (subgroups of)
elementarily free groups
and $\G\subset G_1\times\cdots\times G_n$ a subgroup
such that $L_i:=\G\cap G_i\neq\{1\}$ for $i=1,\dots,n$.

If $L_i$ is finitely generated for $1\le i\le r$
and not finitely generated for $i>r$, then 
there is a subgroup of finite index $\G_0\subset\G$
such that $\G_0 = \G_1\times \G_2$, where $\G_1$
is a direct product of $r$ finitely generated subgroups of
elementarily free groups, and
$H_{k}(\G_2,\Q)$ is infinite dimensional for some $k\le n-r$
(unless $r=n$).
\end{theorem}

In \cite{BHMS} we  deduced a special 
case of Theorem \ref{t:main} from a corresponding
special case of 
Theorem \ref{intersects}.
The latter was proved using an induction  that involved 
the Lyndon-Hochschild-Serre spectral sequence and an analysis of boundary
maps using the Fox calculus. The
 proof of Theorem \ref{t:main} presented here  lays bare 
 some of the
geometry obscured by the algebraic machinery in \cite{BHMS}. 
In
particular,  in the case where
$S$ is not of type $\FP_n$, we construct an 
{\em{explicit family of 
topological cycles}} exhibiting this lack of finiteness.

 One of the ingredients in our proof of Theorem \ref{t:main}
is worthy of particular mention.
Marshall Hall's 
theorem \cite{MH} implies that every non-trivial
element of a finitely generated
free group is primitive (i.e.~generates a free factor)
in a subgroup of finite index. Correspondingly, Peter
Scott \cite{Scott}
proved that any non-trivial element $\gamma$ of the 
fundamental group $\Sigma$ of a closed
surface is represented by a 
non-separating simple closed curve in some finite-sheeted
covering of that surface: algebraically,
this means that there is a subgroup
of finite index $\Sigma_0\subset\Sigma$ such that 
$\Sigma_0$ is an HNN extension $S\ast_C$ where $C$ is
infinite cyclic and the stable letter is $\gamma$.
In Section \ref{s:hasl} we prove the following common generalization of
these theorems.

\begin{theorem}\label{Hall} Suppose that $\G$ splits
as an amalgamated free product $G_1\ast_A G_2$ or
HNN extension $G_1\ast_A$, where $A$ is closed in the
profinite topology on $\G$. If $\gamma\in\Gamma$ is not conjugate
to an element in $G_1$ or $G_2$, then there is a 
subgroup of finite index $\Gamma_0\subset\G$ that splits as
an HNN extension $B\ast_{A'}$, where $A'=A\cap\G_0$
and the stable letter is $\gamma$.
\end{theorem}

This article is organised as follows. In Section \ref{s:lg} we describe those
elements of the structure theory of limit groups that we 
shall need. 
In Section \ref{s:hasl} we prove Theorem \ref{Hall};
our proof relies on the
results in Bass-Serre theory presented in Section \ref{s:bst}. 
In Section \ref{s:dcl} we prove the Double Coset Lemma (Theorem \ref{dc}),
which provides a criterion for homological finiteness in subdirect products
of HNN-extensions. 
In Section \ref{s:general} we marshal the results
of previous sections to prove Theorem \ref{t:main}. 

\section{Limit Groups}\label{s:lg}

Limit groups arise naturally from several points of view. Most
geometrically, they are the finitely generated groups whose Cayley graph can
be obtained as the pointed Gromov-Hausdorff limit of a sequence of 
Cayley graphs of a fixed free group (with a varying choice of generating
set of fixed finite cardinality)
\cite{CG}. Limit groups are precisely those finitely generated
groups $L$ that are {\em{fully residually free}}: 
 for any finite subset $T\subset L$ there exists a homomorphism
from $L$ to a free group that is injective on $T$.
It is in this guise that limit groups were studied extensively by 
Kharlampovich and Myasnikov \cite{KM1}, \cite{KM2}.

The name {\em{limit
group}} was introduced by Zlil Sela to emphasise the fact that these are
precisely the groups that arise when one takes limits of 
stable sequences of homomorphisms $\phi_n:G\to F$, where $G$ is an arbitrary
finitely generated group and $F$ is a free group;  {\em{stable}} means that
for each $g\in G$ either $I_g=\{n\in\N : \phi_n(g)=1\}$ or 
 $J_g=\{n\in\N : \phi_n(g)\neq 1\}$ is finite, and the {\em{limit}} of $(\phi_n)$
is the quotient of $G$ by $\{g\mid |I_g|=\infty\}$.
A good background reference source for limit groups is \cite{BF}.

\subsection{$\omega$-residually free towers}

Our results rely heavily on Sela's version  (\cite{S2}, 1.12; cf
\cite{KM2}) of 
the fundamental theorem  that characterizes limit groups and 
elementarily free
groups in terms of $\omega$-residually free towers. 

\begin{definition}\label{d:sela} 
An $\omega$-rft space of {\em{height}} $h\in\N$
is  defined by induction on $h$. An $\omega$-rft {\em{group}}
is the fundamental group of an $\omega$-rft space.

An $\omega$-rft space of height $0$ is the wedge
(1-point union) of a finite collection of circles, closed hyperbolic
surfaces and tori $\T^n$ (of arbitrary dimension), except that the closed
surface of Euler characteristic $-1$ is excluded\footnote{dropping
this exclusion, and the corresponding one in quadratic blocks,
would not affect the results of our paper}.

An $\omega$-rft space $Y_n$ of height $h$ is obtained from 
an $\omega$-rft space $Y_{h-1}$ of  height $h-1$ by attaching a
{\em{block}} of one of the following types:

\noindent{\em{Abelian:}} $Y_h$ is obtained from $Y_{h-1}\sqcup\T^m$
by identifying a coordinate circle in $\T^m$ with any loop $c$ in 
$Y_{h-1}$ such that $\langle c\rangle\cong\Z$ is a maximal abelian subgroup
of $\pi_1Y_{h-1}$.

\noindent{\em{Quadratic:}} One takes a connected, compact surface $\Sigma$
 that is either a punctured torus or has Euler characteristic at most $-2$, and
obtains $Y_h$  from $Y_{h-1}\sqcup\Sigma$ by identifying each boundary
component of $\Sigma$ with a homotopically non-trivial loop in $Y_{h-1}$; these
identifications must be chosen so that there exists a  retraction 
$r:Y_h\to Y_{h-1}$, and $r_*(\Sigma)\subset\pi_1Y_{h-1}$ must be non-abelian.

\medskip

An $\omega$-rft space is called {\em{hyperbolic}} if no tori are used in its
construction. 
\end{definition}

\begin{theorem}\label{elfree}{\rm (\cite{S6}; see also \cite{KM1},
\cite{KM2})}  A group is elementarily free if and only
if it is the fundamental group of a hyperbolic 
$\o$-rft space. 
\end{theorem}

\begin{thm}\label{embed}{\rm (\cite[1.12]{S2}; see also \cite{KM2})} Limit groups are precisely the
finitely generated subgroups of $\o$-rft groups.
\end{thm}

A useful sketch proof of the latter theorem can be found in \cite{AB}.

This powerful theorem allows one to prove many
interesting facts about limit groups by induction on  height.

\begin{definition}
The {\em{height}} of a limit group $S$ is the minimal height of an $\o$-rft group
that has a subgroup isomorphic to $S$.
\end{definition}

Our approach to understanding an arbitrary limit group $S$
 will be to embed it in the fundamental group of
an $\o$-rft group $L$, take the splitting of $L$
given by Lemma \ref{l:BSbasic} below, and then
decompose $S$ as a graph of groups corresponding to its action
on the Bass-Serre tree of this splitting (subsection \ref{BSss} below).
We refer the reader to \cite{Serre} and \cite{SW} for background on the
Bass-Serre theory of groups acting on trees, which is used extensively
throughout this paper.

\subsection{Elementarily-free  versus hyperbolic limit groups}

A straightforward application of the local gluing lemma 
(\cite{BH}, II.11.3) allows one to deduce from the inductive description given
above that every hyperbolic $\o$-rft supports a locally CAT$(-1)$ metric
and every $\o$-rft supports a locally CAT$(0)$ metric. 
In particular elementarily free groups are hyperbolic. Moreover, an induction
on height, combined with some elementary Bass-Serre theory,
 shows that any
finitely generated subgroup of an elementarily free group $L$ is the fundamental
group of a locally CAT$(-1)$ subcomplex of a  suitable covering
space of the tower space for $L$. Thus:

\begin{lemma} All finitely generated subgroups of elementarily free groups are
hyperbolic limit groups.
\end{lemma}

\begin{remark}
We emphasize that the embedding in 
Theorem \ref{embed} does not
preserve hyperbolicity. For this reason, our proof
of Theorem \ref{t:main} does not extend to arbitrary
hyperbolic limit groups.
\end{remark}

\subsection{Graph of groups decompositions  of elementarily free groups}\label{BSss}

Recall that a graph of groups decomposition  is said to be {\em{$n$-acylindrical}} if 
in the action on the associated Bass-Serre tree, the stabilizer of 
any edge-path of length greater than $n$ is trivial.

\begin{lemma}\label{l:BSbasic}
 If $L$ is the fundamental group of an $\o$-rft  space $Y=Y_h$ of height $h\ge 1$,
then $L$ is the fundamental group of a 2-vertex graph of groups: one
of the vertices is $\pi_1Y_{h-1}$ and the other is a free or free-abelian group
of finite rank at least 2; the edge groups are maximal infinite cyclic subgroups
of the second vertex group. This decomposition is 2-acylindrical.
\end{lemma}

\begin{proof} The splitting described is that which the Seifert-van Kampen
Theorem associates to the decomposition of $Y$ into $Y_{h-1}$ and the
final block in the construction.

If the final block is quadratic, then
the edge groups are precisely the peripheral subgroups of the surface $\Sigma$.
As such they are maximal and form a malnormal family: if $C,C'\subset\Sigma$ are two
cyclic subgroups in the conjugacy classes of two (not necessarily
distinct) edge groups and $x\in\Sigma\smallsetminus
C$, then $x^{-1}Cx\cap C'=\{1\}$.

Any edge-path in the Bass-Serre tree of length greater than
2 must contain a vertex with stabilizer a conjugate of $\Sigma$, and the
intersection of the stabilisers of the incident edges will be of the form
$x^{-1}Cx\cap C'$. Thus the splitting is 2-acylindrical.

Suppose now that the final block in the construction of $Y=Y_h$
is abelian. A straightforward use of the $\o$-residually free
condition yields the following facts about limit groups: first,
if two non-trivial elements 
of an $\o$-residually free group commute and are conjugate,
 then they are equal; second, if $x,y,z$ are non-trivial and 
 $[x,y]=[y,z]=1$, then $[x,z]=1$. It follows from the
 first of these facts that if $A$ is an abelian subgroup
 and $y\in A\cap zAz^{-1}$, then $[y,z]=1$. If $y\neq 1$, it then
 follows from the second fact that $z$ commutes with $A$.
 Thus maximal abelian subgroups of limit groups are malnormal.

By construction, the edge stabilizer in our splitting is 
maximal-abelian in $Y_{h-1}$. Hence it is malnormal in $Y_{h-1}$
and the
splitting is 2-acylindrical.
\end{proof}

\begin{cor}\label{BSsela} If $\G$ is a non-cyclic, finitely generated,
freely indecomposable subgroup of
an elementarily free group $G$, then 
$\G$ is the fundamental group of a 2-acylindrical graph of groups, in which
one of the vertex groups is 
the fundamental group of a compact
surface $\Sigma$ and the incident edge groups are distinct peripheral subgroups 
of $\Sigma$.
\end{cor}
 
\begin{proof}
By Sela's Theorem \ref{elfree}, $G$ is the fundamental group of some $\o$-rft space
$Y$.  If $Y$ can be chosen of height $0$, then $\G$ is the fundamental group
of a closed surface $\Phi$ of Euler characteristic at most $-2$.  In this case the
splitting of $\Phi$ along any nontrivial $2$-sided simple closed curve
induces the desired decomposition of $\G$.  Otherwise, apply Lemma
\ref{l:BSbasic} to $G$, let $T$ be the minimal $\G$-invariant subtree
of the resulting Bass-Serre tree for $G$, and take the resulting graph-of-groups
decomposition for $\G$ with underlying graph $T/\G$.
\end{proof}

\subsection{Subgroup Separability}
Let $\P$ be a class of groups, e.g. free or finite.
Let $\G$ be a group. A subgroup $S\subset\G$
is {\em closed in the pro-$\P$ topology} if  for every $x\not\in S$
there exists a homomorphism
$f:\G\to F$ where $F\in\P$  and $f(x)\not\in f(S)$.

If $\{1\}$ is closed in the pro-$\P$ topology then $\G$
is said to be {\em residually $\P$}.

\noindent{\rm{Notation:}} If all infinite cyclic subgroups $S\subset\G$
are closed in the profinite topology, then we say
that $\G$ is $\Z$-separable.

\begin{prop} In a finitely generated free group $F$, every finitely
generated subgroup is closed in the pro-finite topology.
\end{prop}

\begin{proof} This follows easily from
 Marshall Hall's theorem \cite{MH}, which states that every
finitely generated subgroup is a free factor of a subgroup 
of finite index in $F$.
\end{proof}

\begin{cor} If $S\subset\G$ is closed in the pro-free topology then
it is closed in the pro-finite topology.
\end{cor} 
 
\begin{lemma} If $\G$ is a residually free group with no non-cyclic
abelian subgroups, then $\G$ is $\mathbb Z$-separable.
\end{lemma} 

\begin{proof} 
Let $C=\langle c\rangle\subset\Gamma$  and suppose that $\gamma\in\Gamma$
with $\gamma\not\in C$.
There are two cases to consider: either $[\gamma,c]\neq 1$,
or $H=\langle c,\gamma\rangle$ is cyclic.

In the first case, since $\G$ is residually finite, there exists a homomorphism
$f: \G\to Q$ where $Q$ is finite and $f([\gamma,c])\neq 1$; in particular
$\gamma\not\in f(C)$.

In the second case,  since $\G$ is residually free, there exists a homomorphism
$\phi: \G\to F$ where $F$ is free such that $\phi(\gamma)\neq 1$. Since $\phi|_H$
is injective, $\phi(\gamma)\not\in \phi(C)$, and since $F$ is $\Z$-separable,
there exists a finite quotient such that the image of $\gamma$ does not
lie in the image of $C$.
\end{proof}

\begin{cor} \label{Zsep}
Hyperbolic limit groups (in particular elementarily free groups) are $\Z$-separable.
\end{cor}

A slight variation on the preceding proof
shows that if a group $\G$ is residually free, then
each of its maximal abelian subgroups is closed in the pro-free
topology. Conversations with Zlil Sela and Henry Wilton convinced us that, in the
light of Theorem \ref{embed}, it is not difficult to show that all abelian
subgroups of limit groups are closed in the pro-finite (indeed pro-free)
topology.

\section{Some Bass-Serre Theory}\label{s:bst}

As noted in the previous section, an elementarily free group
acts $2$-acylindrically on the Bass-Serre tree arising from its
$\o$-rf tower decomposition.
In this section we deduce from this that every nontrivial normal subgroup
contains an element that acts hyperbolically on the tree -- a fact
that will be important later.  

Recall that an automorphism $\g$
of a tree $T$ is {\em hyperbolic} if it has no fixed points, and 
that the unique $\g$-minimal subtree of $T$ is then an isometrically
embedded line $A_\g$ called the {\em axis} of $\g$.  In contrast, an
automorphism with at least one fixed point is {\em elliptic}.

\begin{prop}\label{BSp} Let $\G$ be a group acting on a tree $T$.
\begin{enumerate}
\item\label{dishyp} If $\a,\b\in\G$ are hyperbolic  with disjoint axes
$A_\a$ and $A_\b$, then 
$\a\b$ is  hyperbolic and its axis contains the unique
shortest arc from $A_\a$ to $A_\b$.
\item\label{disell} If $\a,\b\in\G$ are elliptic with
${\rm{Fix}}(\a)\cap {\rm{Fix}}(\b)=\emptyset$ then $\a\b$ is
hyperbolic. 
\item\label{convex} If a finite family of convex subsets in $T$ intersects
pairwise, then the intersection of the entire family is non-empty.
\item\label{fghyp} If $\G$ is finitely generated
 then either $\G$ fixes a point of $T$, or else $\G$ contains
hyperbolic isometries.
\item\label{acylhyp} If  
the action of $\G$ is
$n$-acylindrical for some $n\in\mathbb N$, then either $\G$
fixes a point of $T$, or else $\G$ contains hyperbolic isometries.
\item\label{uaxes} If $\G$ contains hyperbolic elements, then the union of the axes of such elements
is the unique minimal $\G$-invariant subtree of $T$.
\end{enumerate}
\end{prop}

\begin{proof}

\ref{dishyp}: Choose an edge $e$ that lies in the arc joining  
$A_\alpha$ to $A_\beta$.
Let $X,Y$ denote the components of $T\setminus\{e\}$ containing $A_\alpha,A_\beta$
respectively.  Note that $\alpha^{\pm 1}(Y\cup e)\subset X$
while $\beta^{\pm 1}(X\cup e)\subset Y$. Thus $e$ is contained in the geodesic
path from $(\alpha\beta)^{-n}(e)$ to $(\alpha\beta)^{m}(e)$ for all $n,m>0$,
and the result follows.

\ref{disell}: A similar argument applies, replacing $A_\alpha,A_\beta$ by
${\rm{Fix}}(\a),{\rm{Fix}}(\b)$ respectively.

\ref{convex}: An inductive argument reduces us to the case of three convex
sets.  Choose a point in each of the three pairwise intersections, and
consider the geodesic triangle in $T$ with these points as vertices.
Each of our three sets
contains one of the  sides of the triangle.
And since we are in a tree, the three sides
have a common point.

\ref{fghyp}: By \ref{disell}, if $\G$ is 
generated by a finite set of elliptics then
either the product of some pair of these generators is
 hyperbolic or else the fixed-point
sets of each pair intersect non-trivially, in which case 
it follows from \ref{convex} that $\G$ has
a fixed point.

\ref{acylhyp}: Without loss of generality, we may assume that  
$T$ is a minimal $\G$-tree. We will
reach a contradiction by assuming that $\Gamma$ has no hyperbolic
elements and no fixed point.

Consider the union $U$  of the sets ${\rm{Fix}}(\g)$ as $\g$ ranges
over $\G\setminus\{1\}$.  We claim that $U$ is convex, hence a subtree.
 Indeed, given $x_1,x_2\in U$ there exist $\g_1,\g_2\in N$ such
that $x_i\in {\rm{Fix}}(\g_i)$ for $i=1,2$, and since $\g_1\g_2$
is not hyperbolic,  ${\rm{Fix}}(\g_1)\cap{\rm{Fix}}(\g_2)\neq\emptyset$
by \ref{disell}.
Thus $[x_1,x_2]$ lies in the convex set  ${\rm{Fix}}(\g_1)\cup {\rm{Fix}}(\g_2)$.
The tree $U$ is nonempty and $\G$-invariant, so $U=T$.

Given  an edge $e$ of $T$, we may choose vertices
$x,y\in T$ in distinct components of $T\ssm\{e\}$ a distance
at least $n$ from $e$. Such vertices exist because  $T$ does 
not have vertices of valence 1 (since removing the edges incident at
such vertices would yield a proper invariant subtree), and therefore
each edge in the complement of   $e$ can be extended  
to an infinite ray in its component of $T\ssm\{e\}$.

Let $\a,\b\in \G$ be non-trivial
elements such that $x\in {\rm{Fix}}(\a)$ and $y\in {\rm{Fix}}(\b)$. Then
 $e$ lies in the convex set ${\rm{Fix}}(\a)\cup {\rm{Fix}}(\b)$. If $e$
lies in the convex set ${\rm{Fix}}(\a)$, then so does 
the segment joining
$x$ to $e$.  But then this segment is fixed by $\alpha$, contradicting the
fact that $\G$ acts $n$-acylindrically.

\ref{uaxes}: This follows from \ref{dishyp}.
\end{proof}

\begin{cor} \label{Nacyc}
If  an  action of a group $\G$ on a tree $T$ is minimal and
$n$-acylindrical for some $n>0$, and $\G$ 
has no global fixed point in $T$,
then every non-trivial normal subgroup 
$N\subset\G$ contains hyperbolic elements.
\end{cor}  

\begin{proof}
Since $N$ is normal in $\G$, the subset ${\rm Fix}(N)$ of $T$ is convex and
$\G$-invariant, so is either the whole of $T$ or is empty. In the first case,
$T$ must have finite diameter by the $n$-acylindrical property.  But then $\G$
has a global fixed point in $T$, contrary to hypothesis.  Hence $N$ has no
global fixed point in $T$, and the result follows from Proposition \ref{BSp}(\ref{acylhyp}). 
\end{proof}
 
\begin{remark} The
 following example illustrates the difficulties
one faces in trying to sharpen the above statement.
Consider the action of the HNN-extension
$B=\langle a,t\mid tat^{-1}=a^2\rangle$
 of the cyclic group $\langle a\rangle$
on its Bass-Serre tree. This action is cocompact with cyclic edge stabilisers.
The normal closure $N$ of $a$ does not contain any hyperbolic
elements.
\end{remark}

\section{Hyperbolics are almost stable letters}\label{s:hasl}

The proof of the following theorem is based on 
an idea introduced by John
 Stallings in the
context of a free group acting freely on its Cayley tree
 \cite{stall2}
(see also \cite[Theorem 2.2]{Scott} and 
\cite[Lemma 15.22]{Hempel}). This result is
a restatement of Theorem \ref{Hall}.

\begin{theorem}\label{MHall}
Let $\G$ be a group that acts 
minimally
on a tree $T$ and let $e$ be an edge whose stabiliser
$A\subset\G$ is closed in the pro-finite topology.

If $t\in \G$ is a 
hyperbolic isometry whose axis contains $e$,
 then there exists a
subgroup $H\subset\G$ of finite index 
such that  $H$ is an HNN-extension with stable letter $t$ and
amalgamated subgroup $A\cap H$.
\end{theorem}

\begin{proof}
Consider the segment
$\l$ of the axis of $t$ that begins with the edge $e$ and ends with
the edge $t(e)$.  Let 
$$e=g_0(e), g_1(e), \dots, g_n(e)=t(e)$$
be the finitely many edges in $\l$ that belong to
  the $\G$-orbit of $e$.  The elements $g_i$ are not
well-defined (unless $A=1$).  But the left cosets $g_iA$ are well-defined
and pairwise distinct.

Since $A\subset\G$ is closed in the profinite topology,
there exists a finite-index normal
subgroup $K\subset\G$ such that the left cosets $g_0(AK),\dots,g_n(AK)$
are pairwise distinct.

Bass-Serre theory expresses $\G$ as the fundamental group
of a graph of groups, where one of the edge-groups is $A$.  Using
the construction described by Scott and Wall in \cite[Proposition 3.6]{SW},
we can construct a classifying space $X$ for $\G$ as a graph of aspherical
spaces.

Thus $\G$ acts freely on the universal cover $\tilde X$ of $X$, and there
is a $\G$-equivariant map $f:\tilde X\to T$ such that $E=f^{-1}(e)$ is
a product (``cylinder") $\tilde X_A\times (0,1)$ for some $K(A,1)$-space $X_A$.
We identify $\tilde X_A$ with $\tilde X_A\times\{\frac12\}\subset E$, and
choose a point $\tilde x_0\in\tilde X_A$ as a base-point for $\tilde X$.
We also choose a path $\tilde\tau$ in $\tilde X$ from $\tilde x_0$ to 
$t(\tilde x_0)$.
Since $\tilde X_A$, and each of its $\G$-translates, has a bicollared neighbourhood in $\tilde X$, we may assume that $\tilde\tau$ is transverse to
$g(\tilde X_A)$ for each $g\in\G$.  Moreover, if we further assume that
the total number of points of intersection of $\tilde\tau$ with
$\bigcup_g~g(\tilde X_A)$ is as small as possible, then $\tilde\tau$
crosses $g(\tilde X_A)$ transversely in a single point for $g=g_1,\dots, g_{n-1}$, meets $g(\tilde X_A)$ in one of its endpoints for $g=g_0=1$
and $g=g_n=t$, and is disjoint from $g(\tilde X_A)$ for all other $g$.

Let $x$ be the image of $\tilde x$ in $X$, and $\tau$ the image of $\tilde\tau$ in $X$. We use $x$ as the base-point for $X$.
 Then $\tau$ is a closed path in $X$ representing
the element $t$ of $\pi_1(X,x)\cong\G$.

Consider the covering space $Y={\tilde X}/K$ of $X$
corresponding to the normal subgroup $K$. Let $E_0$ denote
the image of $E$ in $Y$, and let $X_A'$ be the image
of $\tilde X_A$. The quotient group
$\G/K$ acts on $Y$, and for each $g\in \G$, $E_0$ and $(gK)(E_0)$
either coincide (if $g\in AK$) or are disjoint.  In particular,
$E_0=(g_0K)(E_0),E_1=(g_1K)(E_0),\dots,E_n=(g_nK)(E_0)$ are pairwise disjoint.
Let $y_0$ be the image in $Y$ of $\tilde x_0$, and $\tau'$ the image in
$Y$ of $\tilde\tau$.  Then $\tau'$ is a path from $y_0$ to $(tK)(y_0)$
that intersects each $(g_iK)(X_A')$ in precisely one point.

We cut $E_0$ and $E_n$ along the subspace
$X_A'$  and its translate $(tK)(X_A')=(g_nK)(X_A')$,
creating four (``boundary") 
 subspaces $\partial_-E_0,\partial_-E_n$ and 
$\partial_+E_0, \partial_+E_n$, each
homeomorphic to $X_A'$.  This cutting transforms $Y$ to a space $Y'$,
and $\tau'$ becomes a path in $Y'$ from the copy of $y_0$ in
$\partial_+E_0$ to the copy of
$(tK)(y_0)$ in $\partial_-E_n$.

We form a new space $Z$ from $Y'$ by identifying $\partial_-E_0$ to $\partial_+E_n$,
and $\partial_+E_0$ to $\partial_-E_n$,
using the restriction of the homeomorphism $(tK):E_0\to E_n$.  
Note that this identification in particular identifies
the two endpoints of $\tau'$ to a single point $z$, which we may
regard as the base-point of $Z$.  The image of $\tau'$ in $Z$ is
thus a closed path $\overline\tau$ based at $z$.
The
covering map $Y\to X$ induces a covering map $Z\to X$. (This
$|G:K|$-fold cover is not necessarily connected or regular.)

Consider the component $Z_0$ of $Z$ that contains the copy
$\mathcal A$ of $X_A'$ 
that is the image of $\partial_-E_0$ and $\partial_+E_n$. Then
$z\in {\mathcal A}\subset Z_0$.

The restriction of $Z\to X$
to $Z_0$ is a covering of $X$  corresponding to a
finite index subgroup $H=\pi_1(Z_0,z)\subset\pi_1(X,x)=\G$ {\em{that contains $t$}},
since $\overline\tau$ is a closed path in $Z_0$ which projects to $\tau$.

Since this loop $\overline\tau$
 representing $t\in H=\pi_1(Z_0,z)$ crosses $\mathcal A$ transversely
precisely once,
the Seifert-van Kampen decomposition of $\pi_1(Z_0,z)$ expresses $H$ as 
an HNN-extension with stable letter $t$ and
associated subgroup $\pi_1(\mathcal A)=A\cap H$.
\end{proof}

\begin{cor}\label{mh} Let $\G$ be a group that acts 
minimally
on a tree $T$.  Let $e$ be an edge whose stabiliser
$A\subset\G$ is closed in the pro-finite topology.

If $N\subset\G$ is a normal subgroup that contains a
hyperbolic isometry, then there exists a
subgroup $H\subset\G$ of finite index and an element $t\in N\cap H$
such that  $H$ is an HNN-extension with stable letter $t$ and
amalgamated subgroup $A\cap H$.
\end{cor}

\begin{proof} Since $N$ is normal,
the union $U\subset T$ of the axes of the hyperbolic elements in $N$ is  $\G$-invariant.
But   $U$ is a subtree (Proposition \ref{BSp}(\ref{uaxes})) and  the action of $\G$ is minimal, so
$U=T$.
Hence there exists a hyperbolic element 
$t\in N$ whose axis contains $e$.
\end{proof}

\subsection{The curve-lifting lemma}

We mentioned in the introduction that Theorem \ref{MHall}
 generalizes  Scott's result that every non-trivial
element in the fundamental group of a closed surface can
be represented by a simple closed curve in some finite-sheeted
covering of that surface. In our proof of Theorem \ref{t:main}
we shall need the following refinement of this fact.

\begin{lemma}\label{cl}
Let $\Sigma$ be a compact surface with non-positive Euler characteristic, $X$ a space with $\Z$-separable fundamental group, and $f:\Sigma\to X$ a $\pi_1$-injective map.
If $w$ is a non-trivial element of $\pi_1\Sigma$, then there exists
a finite-sheeted cover $\bar{X}$ of $X$, and a 
simple closed curve
$\alpha$ in the induced cover $\bar\Sigma$ of $\Sigma$, such that
the image of $\alpha$ in $\Sigma$ represents $w$.
\end{lemma}

\begin{proof}
Since $f$ is $\pi_1$-injective, we can identify $\pi_1(\Sigma)$ with
the subgroup $f_*(\pi_1(\Sigma))$ of $\pi_1(X)$.
We do so implicitly throughout the proof without further comment.

Suppose first that $w$ is a proper power of some element $u\in\pi_1\Sigma$;
 say $w=u^n$.  By $\Z$-separability, there is a  subgroup of
 finite
index in $\pi_1(X)$ that contains $w=u^n$ but contains none of  $u,u^2,\dots,u^{n-1}$.  Replacing $X$ and $\Sigma$ by the corresponding finite covers,
we may assume that $w$ is not a proper
power in $\pi_1\Sigma$. We fix a constant-curvature
 Riemannian metric
on $\Sigma$ such that the boundary of $\Sigma$ is totally geodesic.
With respect to this metric, the conjugacy class of $w$ is 
represented by a
closed geodesic $\beta$ with transverse self-intersection.

If $\beta$ is not an embedding, we can express it as a concatenation
of paths $\beta=\beta_1\beta_2\beta_3$ where $\beta_2$ is an 
embedded closed path, which must be
essential in $\Sigma$ since it is a based geodesic in a non-positively
curved metric. Moreover, the conjugacy class represented
by $\beta_2$ has a trivial intersection with
$\langle w\rangle$, because the elements of this cyclic subgroup are represented
by closed geodesics $\beta^n$ which are strictly longer than $\beta_2$.

Since $\pi_1\Sigma\to\pi_1X$ is injective, and $\pi_1X$ is $\Z$-separable,
there is a finite sheeted cover of $X$ such that 
$\beta$ lifts to a closed geodesic in the induced cover of $\Sigma$
 but the corresponding lift of
$\beta_2$ is not closed.  This lift of $\beta$ therefore has fewer
double points than $\beta$.

Repeating this process, we eventually arrive at a finite cover of
$X$ such that $\beta$ lifts to an embedded closed geodesic $\alpha$
in the induced cover of $\Sigma$.
The proof is complete.
\end{proof}

\section{The double-coset lemma}\label{s:dcl}

In the introduction, we explained that our proof of
Theorem \ref{t:main} has the advantage over \cite{BHMS}
that when homological finiteness fails, one can
construct explicit topological cycles that demonstrate
the lack of finiteness. That construction is contained
in the following proof.

\begin{theorem}[Double Coset Lemma]\label{dc}
For $i=1,\dots,n$, let $G_i$ be an HNN-extension with stable letter
$w_i$ and associated subgroups $A_i$ and $B_i$.  Let $G=\prod_{i=1}^n G_i$,
$A=\prod_{i=1}^n A_i$, and let $L\subset G$ be a subgroup containing
$j_i(w_i)$ ($i=1,\dots,n$), where $j_i:G_i\to G$ is the canonical injection.  Suppose also that $p_i(L)=G_i$ for all $i$, where
$p_i:G\to G_i$ is the canonical projection.
Then, the $n$-th homology
group $H_n(L,\Q)$ contains a subgroup isomorphic to
$$\Q\otimes_{\Q L}\Q G\otimes_{\Q A}\Q$$
(the $\Q$ vector space with basis the set of double cosets $\{LgA,g\in G\}$).
\end{theorem}

\begin{proof}
We can construct, for each $i$, a $K(G_i,1)$-complex $X_i$ by starting with
a classifying space for the base of the HNN extension and forming a
 mapping torus corresponding to the given isomorphism $A_i\to B_i$.  
 In this construction,  $w_i$ corresponds to a $1$-cell $W_i$
with both endpoints at the base-point of $X_i$, and $W_i$ appears
only in the boundaries of those $2$-cells corresponding to defining
relations $w_i^{-1}aw_i=b$ for a set of generators $a$ of $A_i$. There
is no loss of generality in assuming that each of the  generators $a,b$ corresponds
to a single 1-cell, and hence the 2-cells involving $w_i$ have attaching maps
of length 4.

Let $X=X_1\times\cdots\times X_n$, let $\tilde{X}$ be the universal
cover of $X$, upon which $G$ acts on the left, and let $\overline{X}=
L\backslash\tilde{X}$ be the covering complex corresponding to
the subgroup $L\subset G$.

We work with cellular chains.

Let $L_i:=\{g\in G_i\mid j_i(g)\in L\}$.  Then $j_i(L_i)$ is the intersection of the kernels of $p_k|_L:L\to G_k$ for $k\ne i$, so
is normal in $L$. Hence $L_i=p_ij_i(L_i)$ is normal in $p_i(L)=G_i$.
Since also $j_i(G_i)$ commutes with $j_k(G_k)$ for $i\ne k$,
it follows that $j_i(L_i)$ is normal in $G$.
But $w_i\in L_i$, so all conjugates of $j_i(w_i)$ belong
to $j_i(L_i)\subset L$.  Hence all the lifts of $W_i$ in $\overline{X}$ are $1$-cycles.

Now let $W$ denote the  cellular $n$-cycle $W_1\times\cdots\times W_n$
 in $X$.
By the above, all the lifts of $W$ to $\overline{X}$ are also
$n$-cycles.  We can identify these lifts with right cosets of $L$
in $G$ as follows.  Choose an $n$-chain $\tilde{W}$ in $\tilde{X}$
that covers $W$. Then the orbit of $W$ under the $G$-action consists of
pairwise distinct $n$-chains $g(\tilde{W})$ for $g\in G$, and 
 $g(\tilde{W})$
and $h(\tilde{W})$ cover the same $n$-chain in $\overline{X}$ if
and only if $Lg=Lh$.

These chains in $\overline{X}$, as has been observed above, are in
fact $n$-cycles, so generate a subgroup $M$
of $H_n(\overline{X},\Q)=H_n(L,\Q)$.
This subgroup is then a homomorphic image of the free $\Q$-module
on the set of right cosets $Lg$ of $L$ in $G$, in other words,
$\Q(L\backslash G)=\Q\otimes_{\Q L}\Q G$.  The kernel of the corresponding
homomorphism is the intersection of this group of $n$-cycles with
the group of (cellular) $(n+1)$-boundaries of $\overline{X}$.

Now the boundary of an $(n+1)$-cell  $\alpha$ of $X$ involves $W$
only if that $(n+1)$-cell is a cube formed as the product of a  $2$-cell
 $\beta$ of some
$X_i$ whose boundary involves $W_i$ and of the $(n-1)$ $1$-cells $W_k$
with $k\ne i$. 
(We arranged in the first paragraph that $\beta$ have an attaching map of
length 4, with two sides corresponding to $W_i$.) Therefore a
lift $\tilde{\alpha}$ of
$\alpha$ in $\tilde{X}$ has the combinatorial structure of a cube, and
 the coefficient of the $n$-dimensional
face $\tilde{W}$ in the cellular $n$-chain
$\partial(\tilde{\alpha})$ is $1-a$ for some $a\in A_i$.

Hence $M$ has a homomorphic image isomorphic to the quotient
of $$\Q\otimes_{\Q L}\Q G$$ by the submodule generated by
 $\{(\Q\otimes_{\Q L}\Q G)(1-a)\mid a\in A\}$.
But this quotient is just $$\Q\otimes_{\Q L}\Q G\otimes_{\Q A}\Q.$$  Since
this is a free $\Q$-module, the epimorphism $$M\to \Q\otimes_{\Q L}\Q G\otimes_{\Q
A}\Q $$
splits, so $M$, and hence also $H_n(L,\Q )$, has a subgroup isomorphic
to $$\Q\otimes_{\Q L}\Q G\otimes_{\Q A}\Q,$$ as claimed.
\end{proof}

We will use the double-coset lemma in the proof of our main result,
where we shall need the following elementary properties of double
cosets in our calculations. 

\begin{lemma}
Let $G$ be a group, $g$ an element of $G$, and $A,B,C$ subgroups of $G$
such that $B\subset C$.  Then the intersection of $C$ and the double coset
$AgB$ is either empty or has the form $(A\cap C)cB$ for some $c\in C$.
\end{lemma}

\begin{proof}
Suppose that $c_1,c_2\in AgB\cap C$.  Then we can write $c_i=a_igb_i$
for $i=1,2$.  Then $(a_2a_1^{-1})=c_2b_2^{-1}b_1c_1^{-1}\in A\cap C$,
so $c_2\in (A\cap C)c_1B$.  Hence, for any $c\in AgB\cap C$, we
have $AgB\cap C\subset(A\cap C)cB$.  The converse inclusion is
immediate, using the equation $AcB=AgB$.
\end{proof}

\begin{cor}\label{c1}
If $A,B,C$ are subgroups of a group $G$ such that $B\subset C$ and
$|A\backslash G/B|<\infty$, then $|(A\cap C)\backslash C/B|<\infty$.
\end{cor}

\begin{lemma}
Let $G,H$ be groups, $A,B$ subgroups of $G$ and $g$ an element of $G$.
If $\phi:G\to H$ is a homomorphism, then $\phi(AgB)=\phi(A)\phi(g)\phi(B)$.
\end{lemma}

\begin{cor}\label{c2}
If $A,B$ are subgroups of a group $G$ such that 
$|A\backslash G/B|<\infty$, and $\phi:G\to H$ is a homomorphism,
then $|\phi(A)\backslash \phi(G)/\phi(B)|<\infty$.
\end{cor}

\section{The Main Theorem}\label{s:general}

In this section we prove the main theorem, Theorem \ref{t:main}, in the
following somewhat stronger form.

\begin{theorem}\label{last}
Let $G_1,\dots,G_n$ be subgroups of
elementarily free groups and let
 $\G\subset G_1\times\cdots\times G_n$ be a subgroup.  Then
the following are equivalent:
\begin{enumerate}
\item[(i)] there exist finitely generated
subgroups $\hat G_i\subset G_i$, for $i=1,\dots,n$,
such that $\G$ is isomorphic to
a finite-index subgroup of $\hat G_1\times\cdots\times\hat G_n$;
\item[(ii)] $\G$ is of type $FP_n(\Q)$;
\item[(iii)] for each $k=1,\dots,n$ and each finite-index subgroup $\G_0\subset\G$,
the $k$-th homology  $H_k(\G_0,\Q)$  is finite dimensional over $\Q$.
\end{enumerate} 
\end{theorem}

\begin{proof}
Finitely generated subgroups of limit groups are limit groups,
and hence of type $FP_\infty$.  Thus any subgroup of finite index
in a direct product of finitely many such groups is also of type
$FP_\infty$, so (i) implies (ii).  It is clear that (ii) implies (iii),
so it remains only to  prove that (iii) implies (i).

Let $G:=G_1\times\cdots\times G_n$, and let $p_i:G\to G_i$ denote
the canonical projection onto the $i$-th factor, for $i=1,\dots,n$.

Replacing each $G_i$ with $p_i(\G)$, we may assume that $p_i(\G)=G_i$
for all $i$.  In particular, each $H_1(G_i,\Q)$ is a homomorphic image
of $H_1(\G,\Q)$, and so finite-dimensinal over $\Q$.  Hence each
$G_i$ is finitely generated, by \cite[Theorem 2]{BHo1}.

By abuse of notation, we identify $G_i$ with the normal
subgroup 
$$\{1\}\times\cdots\times\{1\}\times G_i\times\{1\}\times\cdots\times\{1\}$$
of $G$.  Then we define $L_i=G_i\cap \G$, and note that $L_i$ is normal
in $\G$, and hence that $L_i=p_i(L_i)$ is normal in $G_i=p_i(\G)$.

If $L_n=\{1\}$, then the natural
projection $G_1\times\cdots\times G_n\to G_1\times\cdots\times G_{n-1}$
is injective on $\G$, so $\G$ is isomorphic to a subgroup of 
$G_1\times\cdots\times G_{n-1}\times\{1\}$. An obvious induction
reduces us to the case where all of the
$L_i$ are nontrivial.  We will be done if we can show that
$|G_i:L_i|<\infty$ for each $i$. 

If we replace $G_n$ by a finite-index subgroup $\hat G_n$, and replace each $G_j$
by $\hat G_j=p_jp_n^{-1}(\hat G_n)$ and $\G$ by $\G\cap (\hat G_1\times\dots\times 
\hat G_n)$,
 then neither our
hypotheses nor our desired conclusion is disturbed. (And likewise
with $G_i$ in place of $G_n$.) We shall take advantage of this
freedom several times in 
the sequel without further comment.

We first use this freedom to reduce to the case where none of the $G_i$
is cyclic. (This could also be done by appealing to \cite[Theorem 3]{BHo1}.)
If $G_n$ is cyclic, then $L_n$ has finite index
in $G_n$ since $L_n\ne\{1\}$, so we may assume that $L_n=G_n\cong\Z$.  In this
case $\G$ splits as a direct product $\G'\times\Z$ for some 
$\G'\subset G_1\times\cdots\times G_{n-1}$. The K\"{u}nneth formula
and the homological hypothesis on $\G$ tell us that
 $H_k(\G_0,\Q)$ is finite dimensional
for each $k=1,\dots,n-1$ and for each finite-index subgroup
$\G_0\subset\G'$.  By induction on $n$, we may assume that the theorem
is true for $\G_0$, from which it follows for $\G$.

Henceforth we  assume that none of the $G_i$ is  cyclic.
The structure result for subgroups of 
elementarily free groups, Corollary \ref{BSsela}, tells us that
each $G_i$ is either freely decomposable or can be expressed as the fundamental group
of a $2$-acylindrical graph of groups in which 
one of the vertex groups is
the (nonabelian free)
fundamental group of a surface $\Sigma_i$ with boundary,
and the incident vertex groups are distinct peripheral
subgroups of $\pi_1\Sigma_i$.
  
If $G_i$ is freely decomposable, we regard it as the fundamental group of a nontrivial graph of groups with trivial edge groups (which is
in particular $1$-acylindrical).

The next step concerns those $G_i$ which are freely indecomposable, and
for which $\pi_1\Sigma_i\cap L_i$ is nontrivial.  For such $i$,
we may choose a nontrivial element $a_i$ in $\pi_1\Sigma_i\cap L_i$.
By the curve-lifting lemma, Lemma \ref{cl}, we may assume (after replacing
$G_i$ by a finite-index subgroup), that $a_i$ is represented by a 
simple closed curve in $\Sigma_i$.  Cutting $\Sigma_i$ along this curve
produces a refinement of the graph-of-groups structure of $G_i$,
in which $\<a_i\>$ is an edge group.  Moreover, this refinement remains
$2$-acylindrical.

For those $G_i$ that are freely decomposable, we define $a_i=1$, while for the
remaining $G_i$ we take $a_i$ to be a generator of a peripheral
subgroup of $\pi_1\Sigma_i$ that is an edge
group in the graph-of-groups decomposition of $G_i$.  

In all cases, $G_i$ has a $2$-acylindrical
graph-of-groups decomposition with an edge
group $\<a_i\>$. Moreover $a_i\in L_i$ whenever either $G_i$ is freely decomposable
or $\pi_1\Sigma_i\cap L_i\ne\{1\}$.

By Corollary \ref{Zsep}, $\<a_i\>$ is closed in the profinite topology
on $G_i$; by Corollary \ref{Nacyc}, the normal subgroup $L_i$ contains an
element that acts hyperbolically on the Bass-Serre tree of the decomposition;
hence we may apply Corollary \ref{mh}.  
After replacing the $G_i$
by finite-index subgroups, we may assume that each $G_i$ has an
HNN-decomposition with associated subgroup $\<a_i\>$ and stable letter
$w_i\in L_i$.

By hypothesis $H_n(\G,\Q)$
is finite-dimensional over $\Q$ , so it follows from the double-coset
lemma, Theorem \ref{dc}, that $|\G\backslash G/A|<\infty$, where
$$A=\<a_1\>\times\cdots\times\<a_n\>.$$

We now split the proof into two cases.

\medskip\noindent{\bf Case 1.} Suppose that $a_i\in L_i$ for each $i$.

Since $L_i$ is normal in $G_i$, it follows that $gAg^{-1}\subset \G$
for all $g\in G$, so that $\G gA=\G g$, so $|G:\G|=|\G\backslash G/A|<\infty$,
and the result follows.

\medskip\noindent{\bf Case 2.} Suppose that (possibly after renumbering)
$a_1\notin L_1$.  

Then, by our choice of $a_1$, there is a free
surface group $F=\pi_1\Sigma_1\subset G_1$ with $a_1\in F$ and
$F\cap L_1=\{1\}$.  Define $B=G_1\times\<a_2\>\times\cdots\times\<a_n\>$.
Then $A\subset B$ and $$\left|\G\backslash G/A\right|<\infty,$$ so
Corollary \ref{c1} gives
$$\left|(\G\cap B)\backslash B/A\right|<\infty.$$  Hence Corollary \ref{c2}
gives
$$\left|p_1(\G\cap B)\backslash G_1/\<a_1\> \right| < \infty,$$ since $p_1(B)=G_1$
and $p_1(A)=\<a_1\>$.
But $\<a_1\>\subset F$, so by Corollary \ref{c1} again we have
$$\left|[F\cap p_1(\G\cap B)]\backslash F/\<a_1\>\right|<\infty.$$

Now $G_1$ is normal in $B$ with $B/G_1$ abelian.
Hence $L_1=\G\cap G_1$ is normal in $\G\cap B$ with
$(\G\cap B)/L_1$ abelian.  Hence $L_1=p_1(L_1)$ is normal in
$p_1(\G\cap B)$ with $p_1(\G\cap B)/L_1$ abelian.  Finally,
it follows that $F\cap L_1$ is normal in $F\cap p_1(\G\cap B)$,
with $[F\cap p_1(\G\cap B)]/(F\cap L_1)$ abelian.

But $F\cap L_1=\{1\}$.  Hence $F\cap p_1(\G\cap B)$ is an abelian
subgroup of the free group $F$, and so cyclic: say 
$F\cap p_1(\G\cap B)=\<b\>$.  Then $F$ is a non-abelian free
group, and $a_1,b\in F$ such that $|\<b\>\backslash F/\<a_1\>|<\infty$,
which is absurd.

This contradiction completes the proof.

\end{proof}

\medskip\centerline{\bf Authors' addresses}

\smallskip\begin{center}\begin{tabular}{lll}
Martin R. Bridson &\qquad\qquad & James Howie\\
Department of Mathematics && Department of Mathematics\\
Imperial College London && Heriot-Watt University\\
London SW7 2AZ && Edinburgh EH14 4AS\\
{\tt m.bridson@imperial.ac.uk} && {\tt J.Howie@hw.ac.uk}
\end{tabular}\end{center}

\end{document}